




\input amstex
\documentstyle{amsppt}
\magnification=\magstep1
\NoBlackBoxes


\topmatter

\title Second derivative test for isometric embeddings in $L_p$
\endtitle

\author Alexander Koldobsky \endauthor 
\address Division of Mathematics and Statistics, 
University
of Texas at San Antonio, San Antonio, TX 78249, U.S.A. \endaddress
\email koldobsk\@math.utsa.edu \endemail

\abstract An old problem of P. Levy is to characterize
those Banach spaces which embed isometrically in $L_p.$ 
We show a new criterion in terms of the second
derivative of the norm.
As a consequence we show that, if $M$ is a twice 
differentiable Orlicz function with $M'(0)=M''(0)=0,$ then
the $n$-dimensional Orlicz space $\ell_M^n,\ n\ge 3$ does 
not embed isometrically in $L_p$ with $0<p\le 2.$ 
These results generalize and clear up the recent solution to the 1938 
Schoenberg's problem on positive definite functions.
\endabstract


\rightheadtext{Isometric embeddings in $L_p$}

\thanks Research supported in part 
by the NSF Grant DMS-9531594 \endthanks

\endtopmatter \document \baselineskip=14pt

\head 1. Introduction \endhead

In 1938, Schoenberg \cite{28} asked a question on
positive definite functions an equivalent
formulation of which is as follows: for which 
$p\in (0,2), n\ge 2,\ q>2$ is the space $\ell_q^n$ isometric
to a subspace of $L_p ?$ The solution was completed in 1991
and since then there have appeared a few more proofs (see 
historical remarks below). However, all those solutions were
a little artificial in the sense that the nature of the result
was buried under technical details. The first proof \cite{13}
depended heavily on calculations involving the Fourier transform
of distributions. The later proof of Zastavny \cite{29, 30} was
simpler and led to more general results, but was still quite
technical.

The approach developed in this article not only leads
to further generalizations, but also seems to clear 
things up. To support this ambitious statement, let
us show a very simple argument which, unfortunately,
is false. It was known
to P.Levy \cite {19} that an 
$n$-dimensional normed space $B=(\Bbb R^n,\|\cdot\|)$ 
embeds isometrically in $L_p,\ p>0$ if and only if there
exists a finite Borel measure $\mu$ on the unit sphere
$\Omega$ in $\Bbb R^n$ so that, for every $x\in \Bbb R^n,$ 
$$\|x\|^p =\int_\Omega |(x,\xi)|^p\ d\mu(\xi), \tag{1}$$
where $(x,\xi)$ stands for the scalar product.
Let us formally take the second derivatives by $x_1$ in both
sides of (1). We get
$$p(p-1)\|x\|^{p-2} \big(\|x\|_{x_1}^{'}\big)^2 + 
p\|x\|^{p-1}\|x\|_{x_1^2}^{''} = 
p(p-1) \int_\Omega |(x,\xi)|^{p-2}\ \xi_1^2\ d\mu(\xi). \tag{2}$$
Suppose (and this is the case for the spaces $\ell_q^n$ with
$q>2$) that $\|x\|_{x_1}^{'}(0,x_2,x_3) = 
\|x\|_{x_1^2}^{''}(0,x_2,x_3) = 0$ for every
$(x_2,x_3)\in \Bbb R^2\setminus \{0\}.$ Let us show that under
this condition $B$ can not be isometric to a subspace of $L_p.$
In fact, puting $x_1=0$
in (2) we get zero in the left-hand side. On the other hand, 
the function under the integral in the right-hand side is non-negative,
and the integral can be equal to zero only if the measure $\mu$ is
supported in the hyperplane $\xi_1 = 0.$ But this contradicts to
the assumption that the dimension of $B$ is $n.$

\medbreak

That this reasoning is wrong follows immediately
from the fact that every two-dimensional normed space
(including the spaces $\ell_q^2,\ q>2)$ embeds 
isometrically in $L_p$ for every $p\in (0,1].$
However, if $n\ge 3$ the solution to Schoenberg's problem
is that the spaces $\ell_q^n,\ q>2$ do not embed in 
$L_p$ with $0<p<2,$ so there is a chance that the
argument can be repaired.

\medbreak

The mistake was that for $p\in (0,1]$
the integral in the right-hand side of (2) might diverge
because $p-2 \le -1.$ This, however, can be fixed to a 
certain extent using the so-called 
technique of embedding in $L_{-p}$
introduced in \cite{16}. This technique employs the connection
between the Radon and Fourier transforms to define and study
the representation (1) in the case of negative $p.$ After
relatively short repairments, which clearly distinguish the
two-dimensional case, we get the main result of this article:

\proclaim{Theorem 1} Let $X$ be a three-dimensional normed
space with a normalized basis $e_1, e_2, e_3$ so that: 
\item{(i)} For every fixed
$(x_2,x_3)\in \Bbb R^2\setminus\{0\},$ the 
function $x_1\to \|x_1e_1+x_2e_2+x_3e_3\|$ has continuous 
second derivative everywhere on $\Bbb R,$
and 
$$\|x\|_{x_1}^{'}(0,x_2,x_3) =\|x\|_{x_1^2}^{''}(0,x_2,x_3)= 0;$$ 
\item{(ii)} There exists a constant $K$ so that,
for every $x_1\in \Bbb R$ and every $(x_2,x_3)\in \Bbb R^2$ with
$\|x_2e_2+x_3e_3\|=1,$ we have
$\|x\|_{x_1^2}^{''}(x_1,x_2,x_3) \le K ;$
\item{(iii)} The convergence in the limit 
$\lim_{x_1\to 0} \|x\|_{x_1^2}^{''}
(x_1,x_2,x_3)= 0$ is
uniform with respect to $(x_2,x_3)\in \Bbb R^2$
with $\|x_2e_2+x_3e_3\|=1.$ 

Then, for every $0<p\le 2,$
$X$ is not isometric to a subspace of $L_p.$
\endproclaim

This criterion generalizes the solution to
Schoenberg's problem, 
and it can be applied in some situations where 
previously known criteria do not work. 
In Section 3, we use Theorem 1 to prove 
that Orlicz spaces $\ell_M^n,\ n\ge 3$ do not embed
in $L_p$ with $0<p\le 2$ if $M$ is twice continuously 
differentiable on $[0,\infty)$ and $M'(0)=M''(0)=0.$

\bigbreak

Before we proceed with our repairments, let us 
make a few historical remarks.
The problem of how to check whether a given Banach space
is isometric to a subspace of $L_p$ was raised by P. Levy
\cite {19} in 1937. A well-known fact is that a Banach space embeds
isometrically in a Hilbert space if and only if its norm
satisfies the parallelogram law \cite {6,12}. However, as shown
by Neyman \cite{26}, subspaces of $L_p$ with $p\neq 2$ can not
be characterized by a finite number of equations or inequalities.

P.Levy \cite {19} pointed out that an 
$n$-dimensional normed space $B=(\Bbb R^n,\|\cdot\|)$ 
embeds isometrically in $L_p,\ p>0$ if and only if the
norm admits the representation (1), and,
on the other hand, for $0<p\le 2$ the representation (1)
exists if and only if the function $\exp(-\|x\|^p)$ is 
positive definite, and, hence, is the characteristic function 
of a stable measure. The equivalence of isometric embedding
in $L_p$ with $0\le p \le 2$ and positive definiteness of
the function $\exp(-\|x\|^p)$ was established precisely
by Bretagnolle, Dacunha-Castelle and Krivine \cite {2} who used 
this fact to show that the space $L_q$ embeds isometrically
in $L_p$ if $0<p<q\le 2.$ 

For a long time, the connection with  
stable random vectors and positive definite functions
had been the main source of results on isometric embedding
in $L_p$ (see \cite {1,11,15,17,18,23,24}). However, 
it turns out to be quite difficult to check whether $\exp(-\|x\|^p)$ 
is positive definite for certain norms. For example, the following
1938 Schoenberg's problem \cite{28} was open for more than fifty years:
for which $p\in (0,2)$ is the function  $\exp(-\|x\|_q^p)$ 
positive definite,
where $\|x\|_q$ is the norm of the space $\ell_q^n,\ 2<q\le \infty ?$
As it was mentioned at the beginning of the paper, an equivalent 
formulation asks whether the spaces $\ell_q^n$ embed
in $L_p$ with $0<p\le 2.$ After a few partial results \cite{4}, the 
final answer was given in \cite {22} for $q=\infty$ and in \cite {13}
for $2<q<\infty:$ if $n\ge 3$ the spaces $\ell_q^n,\ 2<q\le \infty$
do not embed isometrically in $L_p,\ 0<p<2.$ A well-known fact
\cite {5,10,20} is that every two-dimensional Banach space embeds in $L_p$
for every $p\in (0,1].$ The spaces $\ell_q^2,\ 2<q\le \infty$ do not
embed in $L_p$ with $1<p\le 2$ (see \cite{4}).

The solution to Schoenberg's problem in \cite{13} was based
on the following Fourier transform characterization of 
subspaces of $L_p$ (see \cite{14}): for every $p>0$ which is 
not an even integer,
an $n$-dimensional  Banach space is isometric to a 
subspace of $L_p$ if and only if the Fourier transform of the 
function $\Gamma(-p/2)\|x\|^p$ is a positive distribution 
on $\Bbb R^n\setminus \{0\}.$ In \cite{3}, this 
criterion was applied to Lorentz spaces. 

Not very long after the paper \cite{13} appeared,
Zastavny \cite {29,30}
proved that a three dimensional space 
is not isometric to a subspace of $L_p$ with $0<p\le 2$ 
if there exists 
a basis $e_1,e_2,e_3$ so that the function
$$(y,z)\mapsto \|xe_1 + ye_2 + ze_3\|^{'}_x (1,y,z)/\|e_1+ye_2+ze_3\|,
\ y,z \in \Bbb R$$
belongs to the space $L_1(\Bbb R^2).$  This criterion provides a
new proof of Schoenberg's conjecture. Zastavny also
showed a stronger result that there are
no non-trivial positive definite functions of the form $f(\|x\|_q).$
(For $q=\infty$ that result was established in \cite{22}; the result for
$2<q<\infty$ was shown independently by Lisitsky \cite{21})

We refer the reader to \cite{15,25} for more 
on isometric embedding of Banach spaces in $L_p$ and its
connections with positive definite functions and isotropic
random vectors. Note that the case $p=1$ is related to the
theory of zonoids. Geometric characterizations of zonoids
and related results of convex geometry can be found in 
\cite{7,9,27}. 

\bigbreak

\head 2. Proof of the second derivative test \endhead

We start with notation and simple remarks.
As usual, we denote by $\Cal S(\Bbb R^n)$ 
the space of rapidly decreasing infinitely differentiable 
functions (test functions) in $\Bbb R^n,$ and 
$\Cal S^{'}(\Bbb R^n)$ is
the space of distributions over $\Cal S(\Bbb R^n).$  
The Fourier transform of a distribution $f\in \Cal S^{'}(\Bbb R^n)$ 
is defined by $\langle \hat{f},\hat{\phi} \rangle = 
(2\pi)^n \langle f,\phi \rangle$ for every
test function $\phi.$  If $p>-1$ and $p$ is not an even integer, then 
the Fourier transform of the function
$h(z)=|z|^p,\ z\in \Bbb R$ is equal to 
$(|z|^p)^{\wedge}(t) = c_p |t|^{-1-p}$ 
(see \cite{8, p.173}), where
$c_p = {{2^{p+1}\sqrt{\pi}\ \Gamma((p+1)/2)}\over{\Gamma(-p/2)}}.$
The well-known connection between the Radon transform and
the Fourier transform is that, for every $\xi\in \Omega,$
the function $t\to \hat\phi(t\xi)$ 
is the Fourier transform of the function 
$z\to R\phi(\xi;z)=\int_{(x, \xi)=z} \phi (x)\, dx$
($R$ stands for the Radon transform). 

\medbreak 
Throughout this section, we remain under the 
conditions and notation of Theorem 1. We denote by
$\|x\|_{x_1}^{'}$ and $\|x\|_{x_1^2}^{''}$  the
first and second partial derivatives by $x_1$ of the norm
$\|x_1e_1+x_2e_2+x_3e_3\|.$ 

\bigbreak

\subheading{Remarks} (i) It is easy to see that, 
for every continuous, homogeneous of degree 
1, positive outside of the origin function $f$ on $\Bbb R^n$ 
and every $\alpha > -n,$ the function $f^{\alpha}$
is locally integrable on $\Bbb R^n.$
In particular, for any $p>0,$ the function
$(x_2,x_3)\to \|x_2e_2+x_3e_3\|^{p-2}$ is locally integrable
on $\Bbb R^2.$
\item{(ii)} A simple consequence of the triangle inequality is that
$-1\le \|x\|_{x_1}^{'} \le 1$ at every point $x\in \Bbb R^3$ with
$(x_2,x_3)\neq 0.$
\item{(iii)} For every fixed $(x_2,x_3)\in \Bbb R^2\setminus \{0\},$ 
$\|x\|$ is a convex differentiable function of $x_1$ whose 
derivative at zero is
equal to zero. Therefore, for every $x=(x_1,x_2,x_3)\in \Bbb R^3,$
we have $\|x\|\ge \|x_2e_2 + x_3e_3\|.$
\item{(iv)} The function $\|x\|_{x_1^2}^{''}$ is non-negative,
homogeneous of degree -1.
Let $K$ be the constant from the condition (ii)
of Theorem 1.
Then, for every $x_1\in \Bbb R$ and every
$(x_2,x_3)\in \Bbb R^2\setminus \{0\},$ the second derivative
$\|x\|_{x_1^2}^{''}(x_1,x_2,x_3)$ is equal to
$${1\over{\|x_2e_2 + x_3e_3\|}}
\|x\|_{x_1^2}^{''}({{x_1}\over{\|x_2e_2 + x_3e_3\|}},
{{x_2}\over{\|x_2e_2 + x_3e_3\|}}, 
{{x_3}\over{\|x_2e_2 + x_3e_3\|}}), $$ 
which is less or equal to $K/\|x_2e_2 + x_3e_3\|$
by the condition (ii) of Theorem 1.

\bigbreak

We are ready to prove the main result of this paper.
\demo{Proof of Theorem 1} For every $0<p_1<p_2\le 2,$ the space 
$L_{p_2}$ embeds isometrically in $L_{p_1},$ so it suffices to prove 
the theorem for $p\in (0,1).$

Suppose that $X$ is isometric to a subspace of $L_p$ with $0<p<1.$
Then, by (1), there exists a measure $\mu$ on the unit sphere $\Omega$ in  
$\Bbb R^3$ so that, for every $x\in \Bbb R^3,$
$$\|x\|^p = \int_\Omega |(x,\xi)|^p\ d\mu(\xi).$$
Applying functions in both sides of the latter equality
to a test function $\phi,$ and using the Fubini theorem and 
the connection between the Radon transform and the Fourier transform,
we get
$$\langle \|x\|^p,\phi \rangle =
\int_{\Bbb R^3} \|x\|^p \phi(x)\ dx =
\int_\Omega d\mu(\xi) \int_{\Bbb R^3} |(x,\xi)|^p \phi(x)\ dx=$$
$$ \int_\Omega d\mu(\xi) \int_{\Bbb R} |t|^p 
\Big(\int_{(x,\xi)=t} \phi(x)\ dx\Big)\ dt =
\int_\Omega \langle |t|^p, R\phi(\xi;t) \rangle\  d\mu(\xi) =$$
$${1\over{(2\pi)^3}}
c_p \int_\Omega \langle |t|^{-1-p}, \hat{\phi}(t\xi) 
\rangle\  d\mu(\xi).\tag{3}$$
By the connection between the Fourier transform and differentiation,
we have $(\partial^2\phi / \partial x_1^2)^{\wedge}(\xi) = -\xi_1^2 \hat{\phi}(\xi),$
and it follows from (3) that
$$\langle \big(\|x\|^p \big)_{x_1^2}^{''},\phi \rangle = 
\langle \|x\|^p, {{\partial^2\phi}\over{\partial x_1^2}} \rangle =
-{1\over{(2\pi)^3}}c_p \int_\Omega \xi_1^2\ d\mu(\xi) 
\int_{\Bbb R} |t|^{1-p}\hat{\phi}(t\xi)\ dt.\tag{4}$$

Let $\phi_n(x_1,x_2,x_3)= h_n(x_1)u(x_2,x_3),$ where
$$h_n(x_1) = (n/\sqrt{2\pi})\exp(-x_1^2n^2/2),
\ u(x_2,x_3) = (1/(2\pi)) \exp(-(x_2^2+x_3^2)/2).$$
Then $\hat{\phi_n}(\xi_1,\xi_2,\xi_3)=
\exp(-{{\xi_1^2}\over {2n^2}}) \exp(-(\xi_2^2+\xi_3^2)/2).$
Clearly, $\int_{\Bbb R} h_n(x_1)\ dx_1 = 1,$ and
$\lim_{n\to \infty} \int_{|x_1|>\delta} h_n(x_1)\ dx_1 =0$
for every $\delta>0.$

Note that the number $c_p$ is negative for every $p\in (0,2)$ 
(see the expression for $c_p$ in the beginning of this
section). 
Applying (4) to the function $\phi_n,$ we get 
$$\langle \big(\|x\|^p \big)_{x_1^2}^{''},\phi_n \rangle =
-2^{-p\over 2} \Gamma(1-{p\over 2}){1\over{(2\pi)^3}} c_p 
\int_\Omega \xi_1^2({{\xi_1^2}\over{n^2}}+ 
\xi_2^2+\xi_3^2)^{{{p-2}\over 2}}
\ d\mu(\xi) \ge $$
$$-2^{-{p\over 2}} \Gamma(1-{p\over 2}) {1\over{(2\pi)^3}}
c_p \int_\Omega \xi_1^2\ d\mu(\xi).\tag{5}$$
By Lemma 1 below, we can make the left-hand side of (5)
as small as we want.
Therefore, $\int_\Omega \xi_1^2\ d\mu(\xi)=0,$ and the measure 
$\mu$ is supported in the hyperplane $\xi_1=0.$
This contradicts 
to the assumption that the space
$X$ is three-dimensional, so $X$ is not isometric
to a subspace of $L_p.$ \qed \enddemo
 
\bigbreak

In the following lemma the norm has the same properties as in 
Theorem 1 and the functions $\phi_n,\ h_n,\ u$ are the same as 
in the proof of Theorem 1.

\proclaim{Lemma 1} For every $\epsilon >0$ there exists $N\in \Bbb N$ 
so that, for every $n>N,$ $ \langle \big(\|x\|^p \big)_{x_1^2}^{''}\ ,
\phi_n \rangle \le \epsilon.$\endproclaim

\demo{Proof} First we use the Fubini theorem to show that
$$\langle \big(\|x\|^p \big)_{x_1^2}^{''}\ ,\phi_n \rangle =
\langle \|x\|^p, 
{{\partial^2\phi_n}\over{\partial x_1^2}} \rangle =$$
$$\int_{\Bbb R^3} \|x\|^p 
{{\partial^2\phi_n}\over{\partial x_1^2}}\ dx =
\int_{\Bbb R^2\setminus\{0\}} u(x_2,x_3)\ dx_2\ dx_3
\Big( \int_{\Bbb R} \|x\|^p 
{{\partial^2 h_n}\over{\partial x_1^2}}\ dx_1 \Big) =$$
$$ \int_{\Bbb R^2\setminus\{0\}} 
\langle \big(\|x\|^p \big)_{x_1^2}^{''}\ ,h_n \rangle
u(x_2,x_3)\ dx_2\ dx_3.\tag{6}$$

We restrict the outer integral to $\Bbb R^2\setminus\{0\}$
(which does not change the value of the integral)
to formally exclude the point $(x_2,x_3)=0,$ where the derivative 
$\big(\|x\|^p \big)_{x_1^2}^{''}$ is a distribution
in terms of the $\delta$-function of $x_1.$

By the condition (i) of Theorem 1 and Remarks (i), (ii), (iv),
for every fixed $(x_2,x_3)\in \Bbb R^2\setminus\{0\},$ 
the derivative $\big(\|x\|^p \big)_{x_1^2}^{''}$ is
a locally integrable continuous function of the variable $x_1$
on $\Bbb R,$ so the expression in (6) can be written as 
$$ \int_{\Bbb R^2\setminus\{0\}} 
\Big(\int_{\Bbb R} \big(\|x\|^p \big)_{x_1^2}^{''} h_n(x_1)\ dx_1
\Big) u(x_2,x_3)\ dx_2\ dx_3 = \tag{7}$$
$$\int_{\Bbb R^2\setminus\{0\}} \Big(\int_{\Bbb R}
\big(p(p-1)\|x\|^{p-2} 
\big(\|x\|_{x_1}^{'}\big)^2 + 
p\|x\|^{p-1}\|x\|_{x_1^2}^{''}\big)
h_n(x_1)\ dx_1 \Big) u(x_2,x_3)\ dx_2 dx_3.$$
Since $p<1,$ the first term under the intergal by $x_1$
is negative, so, to prove Lemma 1, it suffices to show that 
$$\lim_{n\to \infty} \int_{\Bbb R\times(\Bbb R^2\setminus\{0\})} 
\|x\|^{p-1}\|x\|_{x_1^2}^{''}
\ h_n(x_1)\  u(x_2,x_3)\ dx = 0.$$

Let $\epsilon > 0.$ By Remark (i), the function
$\|x_2e_2+x_3e_3\|^{p-2}$ is locally integrable in $\Bbb R^2,$
so 
$$L= \int_{\Bbb R^2}\|x_2e_2+x_3e_3\|^{p-2} u(x_2,x_3)
\ dx_2\ dx_3 < \infty.$$
Besides, there exists  $c>0$ so that
$$\int_{\{(x_2,x_3): \|x_2e_2+x_3e_3\|<c\}}
\|x_2e_2+x_3e_3\|^{p-2} u(x_2,x_3)
\ dx_2 \ dx_3 <  {{\epsilon}\over{3K}},\tag{8}$$
where $K$ is the number defined in Remark (iv). By Remark (iv),
we have
$$\int_{\Bbb R \times \{(x_2,x_3): \|x_2e_2+x_3e_3\|<c\}}
\|x\|^{p-1} \|x\|_{x_1^2}^{''}\ h_n(x_1)\ u(x_2,x_3)\ dx \le$$
$$K\int_{\Bbb R\times \{(x_2,x_3): \|x_2e_2+x_3e_3\|<c\}} 
\|x\|^{p-1} \|x_2e_2+x_3e_3\|^{-1} 
h_n(x_1) u(x_2,x_3)\ dx < {\epsilon\over 3},\tag{9}$$
where we apply Remark (iii) to show that
$\|x\|^{p-1}\le \|x_2e_2+x_3e_3\|^{p-1}$
(note that $p-1<0),$ and then use
the inequality (8) and the fact that $\int_{\Bbb R} h_n(t)\ dt = 1.$

By the condition (iii) of Theorem 1, $\lim_{x_1\to 0} \|x\|_{x_1^2}^{''}
(x_1,x_2,x_3) = 0$ uniformly with respect to $(x_2,x_3)\in \Bbb R^2$
with $\|x_2e_2+x_3e_3\|=1,$  so there exists $\delta>0$ so that if
$|x_1|<\delta$ then  $\|x\|_{x_1^2}^{''}(x_1,x_2,x_3) < 
{\epsilon\over {3L}}$
for every $(x_2,x_3)\in \Bbb R^2$ with $\|x_2e_2+x_3e_3\|=1.$
The second derivative of the norm is a homogeneous function
of degree -1. Therefore, if $(x_2,x_3)\in \Bbb R^2\setminus\{0\},$
and ${{|x_1|}\over{\|x_2e_2 + x_3e_3\|}} < \delta,$ then
$$\|x\|_{x_1^2}^{''}(x_1,x_2,x_3) =
{1\over{\|x_2e_2+x_3e_3\|}} \|x\|_{x_1^2}^{''}
({x\over{\|x_2e_2+x_3e_3\|}})
< {\epsilon\over{3L\|x_2e_2 + x_3e_3\|}}.
\tag{10}$$

Denote by  $A_1, A_2$ the sets 
${\{x_1: |x_1|<\delta c\} \times \{(x_2,x_3): \|x_2e_2+x_3e_3\|\ge c\}}$
and ${\{x_1: |x_1|>\delta c\} \times \{(x_2,x_3): \|x_2e_2+x_3e_3\|\ge c\}}$
in $\Bbb R^3.$ Then, by (10), Remark (iii) and the fact that
$\int_{\Bbb R} h_n(t)\ dt = 1,$
$$\int_{A_1} \|x\|^{p-1} \|x\|_{x_1^2}^{''} h_n(x_1)\ u(x_2,x_3)\ dx \le$$
$${\epsilon\over {3L}}
\int_{A_1} \|x_2e_2+x_3e_3\|^{p-2}h_n(x_1)\ u(x_2,x_3)\ dx \le 
{\epsilon\over 3}. \tag{11}$$

Finally, suppose that $n$ is large enough so that 
$\int_{|x_1|>\delta c} h_n(x_1)\ dx_1 < \epsilon/(3KL).$
Then we use the estimate of Remark (iv) and Remark (iii) to show that
$$\int_{A_2} \|x\|^{p-1} \|x\|_{x_1^2}^{''} h_n(x_1)\ u(x_2,x_3)\ dx \le$$
$$K \int_{A_2} \|x_2e_2+x_3e_3\|^{p-2}h_n(x_1)\ u(x_2,x_3)\ dx \le
{\epsilon\over 3}.\tag{12}$$

Since $\epsilon$ is an arbitrary positive number, 
the result of Lemma 1 follows from (9), (11), (12)
and the fact that the second derivative of 
a convex function is non-negative.
\qed \enddemo

\bigbreak

\subheading{Remark} Note that the statement of Theorem 1 
is not true for two-dimensional spaces all of which embed
isometrically in $L_p$ for every $p\in (0,1].$ The reason 
that the proof of Theorem 1 does not work for two-dimensional
spaces is that, unlike the function $\|x_2e_2+x_3e_3\|^{p-2}$ 
on $\Bbb R^2,$ the function $|x_2|^{p-2}$ is not locally 
integrable on $\Bbb R$ if $0<p\le 1.$

\head 3. Application to Orlicz spaces \endhead

An Orlicz function $M$ is a non-decreasing convex function
on $[0,\infty)$ such that $M(0)=0$ and 
$M(t)>0$ for every $t>0.$

For an Orlicz function $M,$ the norm of the $n$-dimensional
Orlicz space $\ell_M^n$ is defined by the equality
$\sum_{k=1}^n M(|x_k|/\|x\|) = 1,\ x\in \Bbb R^n\setminus \{0\}.$

\proclaim{Theorem 2} Let $M$ be an Orlicz function so that
$M\in C^2([0,\infty)),\ M'(0) = M''(0) =0.$ Then, for every
$0<p\le 2,$ the three-dimensional Orlicz space $\ell_M^3$
does not embed isometrically in $L_p.$ \endproclaim

\demo{Proof} We are going to show that the norm 
of the space $\ell_M^3$ satisfies the conditions of Theorem 1.
Since the Orlicz norm is an even function with respect 
to each variable, it suffices to consider
the points $x=(x_1,x_2,x_3)$ with non-negative coordinates. 
We denote by $e_1,e_2,e_3$ the standard normalized basis
in $\ell_M^3.$

The function $M'$ is non-decreasing, continuous on 
$[0,\infty)$ and $M'(0)=0.$ Since $M(0)=0$ and $M(t)>0$ for every
$t>0,$ the function $M'$ can not be equal to zero on
an interval, so $M'(t)>0$ for every $t>0.$

Let $x=(x_1,x_2,x_3)$ with $(x_2,x_3)\neq 0.$ Then
one of the numbers $x_2 M'(x_2/\|x\|)$ or  $x_3 M'(x_3/\|x\|)$
is positive. By implicit differentiation,
$$\|x\|_{x_1}^{'} = {{\|x\| M'({{x_1}\over {\|x\|}})}\over
{x_1 M'({{x_1}\over {\|x\|}}) + x_2 M'({{x_2}\over {\|x\|}}) +
x_3 M'({{x_3}\over {\|x\|}})}}.\tag{13}$$
Also,
$$\|x\|_{x_1^2}^{''} = {{(\|x\|-x_1\|x\|_{x_1}^{'})^2 
M''({{x_1}\over {\|x\|}})+ x_2^2 (\|x\|_{x_1}^{'})^2
M''({{x_2}\over {\|x\|}})+ x_3^2 (\|x\|_{x_1}^{'})^2 
M''({{x_3}\over {\|x\|}})}\over 
{\|x\|^2 \big(x_1 M'({{x_1}\over {\|x\|}}) + 
x_2 M'({{x_2}\over {\|x\|}}) +
x_3 M'({{x_3}\over {\|x\|}})\big)}}.\tag{14}$$
The condition (i) of Theorem 1 follows from the fact
that $M'(0)=M''(0)=0.$

Let us show that the norm satisfies the condition (ii)
of Theorem 1.
Denote by $c=\min\{x_2 M'(x_2/2) +
x_3 M'(x_3/2): \|x_2e_2+x_3e_3\|=1,\ x_2,x_3\ge 0.\}$
Since $M'$ is a continuous function and $M'(t)>0$ for $t>0,$
we have $c>0.$ Let $d=\max_{t\in [0,1]} M''(t).$

Clearly, $x_i\le \|x\|,\ i=1,2,3.$
Therefore, using also Remark (ii) and positivity
of $ \|x\|_{x_1}^{'},$ we have 
$0\le \|x\| - x_1 \|x\|_{x_1}^{'} \le \|x\|,$ and 
$(\|x\| - x_1 \|x\|_{x_1}^{'})^2/\|x\|^2 \le 1.$ 

Consider any $x_2,x_3\ge 0$ with $\|x_2e_2+x_3e_3\|=1.$
Then $x_2,x_3 \le 1.$ 
If $x_1\in [0,1]$ then $1\le \|x\|\le 2,$ hence,
$x_i/\|x\|\ge (x_i/2),\ i=1,2,3,$ 
We get from (14) that $\|x\|_{x_1^2}^{''} \le 3d/c.$

If $x_1 >1$ then $x_1/\|x\|>1/2,$ and (14) implies
$\|x\|_{x_1^2}^{''} \le 3d/M'(1/2),$ because $M'$ is
an increasing function.

Since in both cases $x_1\in [0,1]$ and $x_1>1$ we estimated
the second derivative by constants which do not depend on
the choice of $x_2,x_3$ with $\|x_2e_2+x_3e_3\|=1,$ 
we get the condition (ii) of Theorem 1.

Finally, let us show the condition (iii) of Theorem 1.
Let $c$ be as defined above. Denote by 
Since $M'$ is a continuous increasing function whose 
value at zero is zero,
for every $\epsilon>0$ there exists $\delta>0$ so that
$M'(t)< c\epsilon/2$ if $x_1<\delta.$ Then, by (13),
for every $x_2,x_3$ with $\|x_2e_2+x_3e_3\|=1,$
$\|x\|_{x_1}^{'}(x_1,x_2,x_3) \le \epsilon$ 
if $x_1<\min(1,\delta),$
which proves that the first derivative converges to zero
uniformly. Similarly, we use (14), the uniform 
convergence of the first derivative and the fact that 
$M''$ is continuous and $M''(0)=0$ to prove that the
second derivative of the norm also converges to zero
uniformly.
\qed \enddemo

\bigbreak

Clearly, the class of Orlicz functions satisfying the
conditions of Theorem 2
includes all the functions $M(t)=|t|^q,\ q>2.$ Also, 
if $M$ satisfies the conditions of Theorem 2, the spaces
$\ell_M^n$ with $n\ge 3$ (as well as the infinite dimensional
spaces $\ell_M,\ L_M([0,1]))$ do not embed 
isometrically in $L_p$ with
$0<p\le 2.$ Therefore,
for every $0<p\le 2$ and the norm of each of those spaces, 
the function $\exp(-\|x\|^p)$ is not positive definite.
The latter result generalizes the solution to the Schoenberg's
problem.

\enddocument